\RequirePackage{ifpdf}
\ifpdf 
\documentclass[pdftex]{sigma}
\else
\documentclass{sigma}
\fi

\usepackage{shadow}

\begin{document}

\renewcommand{\PaperNumber}{089}

\def\pr{{}^{\,\prime}}

\def\ve{\varepsilon}
\def\ss{\scriptstyle}
\def\sss{\scriptscriptstyle}
\newcommand{\sy}[1]{(\hspace{-.01mm}\underline{#1})}
\newcommand{\syt}[1]{(\hspace{-.01mm}\underline{#1})^{\!\rm T}}
\newcommand{\sytt}[1]{(\hspace{-.01mm}\underline{#1})^{\!\rm T\!T}}
\newcommand{\syg}[1]{(\hspace{-.01mm}\underline{#1})^{\rm \gamma}}
\newcommand{\sygt}[1]{(\hspace{-.01mm}\underline{#1})^{\rm \gamma T}}

\newcommand{\DB}[1]{\Delta^{\!(#1)}}
\newcommand{\DF}[1]{{\cal D}^{(#1)}}
\newcommand{\lDF}[1]{\overleftarrow{\cal D}^{(#1)}}

\newcommand{\Dp}[1]{D^{(#1)}}
\renewcommand{\S}[1]{S^{(#1/2)}}

\def\TT{{\rm TT}}
\def\N{{\cal N}}

\def\GG{{\cal G\hspace{-2.3mm}G}}
\def\Gg{\mathbb G}
\def\E{{\cal E}}
\def\B{{\cal B}}
\def\lrd{\stackrel{\leftrightarrow}{\partial}}

\def\chib{\overline\chi}
\def\psib{\overline\psi\hspace{-2.6mm}\phantom{\psi}}
\def\k{\kappa}
\def\tid{\hspace{.3mm}}
\def\l{\lambda}
\def\ve{\varepsilon}
\def\D{{\cal D}}
\def\M{{\cal M}}
\def\J{{\cal J}}
\def\r{\rho}
\def\s{\sigma}
\def\t{\tau}
\def\ph#1{\phantom{#1}}
\def\ag#1{\gamma({\scriptstyle #1})}
\def\wh#1{\widehat{#1}}
\def\wt#1{\widetilde{#1}}
\def\mt{{\tilde \mu}}
\def\nt{{\tilde \nu}}
\def\rt{{\tilde \rho}}
\def\st{{\tilde \sigma}}
\def\ol#1{\overline{#1}}
\def\sss{\scriptscriptstyle}
\def\ts{\textstyle}
\def\d{{\bf d}}
\def\m{\mu}
\def\n{\nu}
\def\la{\langle}
\def\ra{\rangle}
\def\e{\epsilon}
\def\scirc{\!{\scriptstyle \circ}}
\def\F{{\cal F}}
\def\R{{\cal R}}
\def\psitb{\overline{\widetilde\psi}}
\def\vphi{\varphi}
\def\psid{\psi^\dagger}

\def\beqa{\begin{equation}\begin{array}{l}}
\def\eeqa{\end{array}\end{equation}}
\def\eqlab#1{\label{eq:#1}}
\def\figlab#1{\label{fig:#1}}
\def\tablab#1{\label{tab:#1}}
\def\seclab#1{\label{sec:#1}}
\def\eqn#1{(\ref{#1})}
\def\eref#1{(\ref{eq:#1})}
\def\eqref#1{eq.~(\ref{eq:#1})}
\def\Eqref#1{Eq.~(\ref{eq:#1})}
\def\figref#1{fig.~\ref{fig:#1}}
\def\Figref#1{Fig.~\ref{fig:#1}}
\def\tabref#1{\ref{tab:#1}}
\def\Tabref#1{Table \ref{tab:#1}}
\def\secref#1{Section \ref{sec:#1}}
\def\sla#1{#1 \hspace{-2.7mm} \slash}
\def\slad{\partial \hspace{-2.2mm} \slash}
\def\slap{p \hspace{-1.8mm} \slash}
\def\boxfrac#1#2{\mbox{\small{$\frac{#1}{#2}$}}}
\def\half{\mbox{\small{$\frac{1}{2}$}}}
\def\sfrac#1{\mbox{\small{$\frac{1}{#1}$}}}
\def\thalf{\mbox{\small{$\frac{3}{2}$}}}
\def\quarter{\mbox{\small{$\frac{1}{4}$}}}
\def\third{\mbox{\small{$\frac{1}{3}$}}}
\def\sixth{\mbox{\small{$\frac{1}{6}$}}}

\def\a{\alpha}
\def\ga{\gamma} \def\G{{\it\Gamma}} \def\g{\gamma}
\def\de{\delta} \def\De{\Delta}
\def\veps{\varepsilon}  \def\eps{\epsilon}
\def\kp{\kappa}
\def\L{{\it\Lambda}}
\def\Pit{{\it\Pi}}
\def\Psit{{\it\Psi}}
\def\si{\sigma} \def\Si{{\it\Sigma}}
\def\th{\theta} \def\vth{\vartheta} \def\Th{\Theta}
\def\w{\omega} \def\W{\Omega} \def\hw{\hat{\omega}}
\def\vfi{\varphi}\def\vphi{\varphi}
\def\z{\zeta}
\def\bra{\langle} \def\ket{\rangle}
\def\dd{{\rm d}}
\def\pa{\partial}
\def\vrho{\varrho}

\def\ie{{i.e., }}
\def\eg{{e.g.\ }}
\def\cl{\centerline}

\def\pa{\partial}
\def\rarr{\rightarrow}
\def\nn{\nonumber}

\def\psibar{\overline{\psi}}
\def\Gbar{\overline{G}}

\def\psidag{\psi^\dagger}
\def\Gdag{G^\dagger}

\def\N{{\bf N}}
\def\TR{{\bf tr}}
\def\G{{\bf g}}
\def\DIV{{\bf div}}
\def\GRAD{{\bf grad}}
\def\ord{{\bf ord}}

\def\Nt{{\cal N}}
\def\Tt{{\cal T}}
\def\c{{\hspace{.3mm}\bf c\hspace{.3mm}}}
\def\Real{\mathbb R}

\newcommand{\comment}[1]{{\bf [#1]}}
\newcommand{\ul}[1]{\underline{#1}}
\newcommand{\bm}[1]{\mbox{\boldmath $ #1 $}}

\FirstPageHeading

\renewcommand{\thefootnote}{$\star$}

\ShortArticleName{The Symmetric Tensor Lichnerowicz Algebra}

\ArticleName{The Symmetric Tensor Lichnerowicz Algebra\\ and a Novel Associative Fourier--Jacobi Algebra\footnote{This paper is a
contribution to the Proceedings of the 2007 Midwest
Geometry Conference in honor of Thomas~P.\ Branson. The full collection is available at
\href{http://www.emis.de/journals/SIGMA/MGC2007.html}{http://www.emis.de/journals/SIGMA/MGC2007.html}}}

\Author{Karl HALLOWELL and Andrew WALDRON}

\AuthorNameForHeading{K. Hallowell and A. Waldron}

\Address{Department of Mathematics, University of California,
            Davis CA 95616, USA}

\Email{\href{hallowell@math.ucdavis.edu}{hallowell@math.ucdavis.edu}, \href{wally@math.ucdavis.edu}{wally@math.ucdavis.edu}}

\ArticleDates{Received July 21, 2007; Published online September 13, 2007}

\Abstract{Lichnerowicz's algebra of dif\/ferential geometric operators acting on
symmetric tensors can be obtained  from generalized geodesic motion of an observer carrying
a complex tangent vector. This relation is based upon quantizing the classical evolution
equations, and identifying wavefunctions with sections of the symmetric tensor bundle
and Noether charges with geometric operators.  In general curved spaces these operators
obey a deformation of the Fourier--Jacobi Lie algebra of $sp(2,{\mathbb R})$. These results have
already been generalized by the authors to arbitrary tensor and spinor bundles
using supersymmetric quantum mechanical models and have also been applied to
the theory of higher spin particles. These Proceedings review these results in their simplest,
symmetric tensor setting. New results on a novel and extremely useful reformulation of the
rank~2 deformation of the Fourier--Jacobi Lie algebra in terms of an associative algebra
are also presented. This new algebra was originally motivated by studies of operator orderings
in enveloping algebras. It provides a new method that is superior in many respects to common
techniques such as Weyl or normal ordering.}

\Keywords{symmetric tensors; Fourier--Jacobi algebras; higher spins; operator orderings}

\Classification{51P05; 53Z05; 53B21; 70H33; 81R99; 81S10}

\renewcommand{\thefootnote}{\arabic{footnote}}
\setcounter{footnote}{0}

\section{Introduction}

The study of geometry using f\/irst quantized particle models has a long history.
Notable examp\-les are the study of Pontryagin classes and Morse theory in
terms of ${\cal N}=1$ and ${\cal N}=2$ supersymmetric quantum mechanical models~\cite{Alvarez-Gaume:1983ig,Witten:1982im}. The supercharges of those models correspond to
Dirac, and exterior derivative and codif\/ferential operators acting on spinors and forms, respectively.
The model we concentrate on here describes gradient and divergence operators
acting on symmetric tensors and therefore involves {\it no} supersymmetries at all.
Hence, even though the symmetries of this model are analogous to supersymmetries,
no knowledge of superalgebras is required to read these Proceedings.
All the above models f\/it into a very general class of orthosymplectic spinning particle theories
studied in detail by the authors in~\cite{Hallowell:2007qk}. Spinors, dif\/feren\-tial forms,
multiforms~\cite{Olver,Henneaux,Senovilla,Bekaert:2002dt,deMedeiros:2003dc}
and symmetric tensors\footnote{See~\cite{Lecomte} for a f\/lat space discussion of
the symmetric tensor theory and~\cite{Labastida:1987kw,Vasiliev:1988xc} for its origins
in higher spin theories.}~\cite{Lichnerowicz:1961,Hallowell}  are all f\/itted into a single framework
in that work. Here we focus on the symmetric tensor case, both for its simplicity,
and because we want to present new results on the symmetric tensor Lichnerowicz algebra
developed in~\cite{Hallowell}.

The underlying classical system is geodesic motion on a Riemannian manifold along with paral\-lel
transport of a complex tangent vector. This is described by a pair of
ordinary dif\/feren\-tial equations to which we add further curvature couplings designed
to maximize the set of constants of the motion. Of particular interest are symmetries interchanging
the vector tangent to the manifold with the tangent vector to the geodesic.
These are analogous to supersymmet\-ries and correspond to gradient and divergence
operators. This correspondence is achieved by quantizing the model. The complex tangent
vector describes spinning degrees of freedom so that wavefunctions are sections
of the symmetric tensor bundle. The Noether charges of the theory become operators
on these sections. In particular, the Hamiltonian is a curvature modif\/ied Laplace operator. In fact, it is precisely the wave operator acting on symmetric tensors
introduced some time ago by Lichnerowicz~\cite{Lichnerowicz:1961} on the basis of
its algebraic properties on symmetric spaces~\cite{Hallowell}. Moreover, the set of all Noether
charges obey a deformation of the Fourier--Jacobi Lie algebra $sp(2,{\mathbb R})^J$.
The classical model is described in Section~\ref{classical}, while its quantization and
relation to geometry are given in Section~\ref{quantum}.

Applications, such as higher spin theories~\cite{Hallowell,Bekaert:2005vh}, call for
expressions in  the  universal enveloping algebra ${\cal U}(sp(2,{\mathbb R})^J)$
involving arbitrarily high powers of the generators. Manipulating these expressions
requires a standard ordering, oft used examples being Weyl ordering
(averaging over operator orderings) or normal ordering (based on a choice of polarization
such that certain operators are moved preferentially to the right, say). In a study of partially massless
higher spins~\cite{Deser:2001pe}, we found a new operator ordering scheme to be particularly
advantageous~\cite{Hallowell}. The key idea is to rewrite generators of the $sp(2,{\mathbb R})$
subalgebra, wherever possible, as po\-wers of Cartan elements or the quadratic Casimir operator.
Immediately, this scheme runs into a~dif\/f\/iculty, namely that the remaining $sp(2,{\mathbb R})^J$
generators do not have a simple commutation relation with the quadratic $sp(2,{\mathbb R})$
Casimir. This problem is solved by a trick: we introduce a~certain square root
of the quadratic Casimir whose {\it r\^ole} is to measure how far states are
from being highest weight. Then we use this square root operator to construct modif\/ied
versions of the remaining $sp(2,{\mathbb R})^J$ generators. Instead of a simple Lie algebra,
we then obtain an elegant associative algeb\-ra, which we denote $\wt{\, \cal U \,}\!(sp(2,{\mathbb R}))$,
with relations allowing elements to be easi\-ly reordered. This algebra is described and derived
in detail in Section~\ref{FJalgebra}. The f\/inal Section discusses applications and our conclusions.

\section{The classical model}
\label{classical}

\begin{figure}[t]
\centerline{\includegraphics[height=6cm]{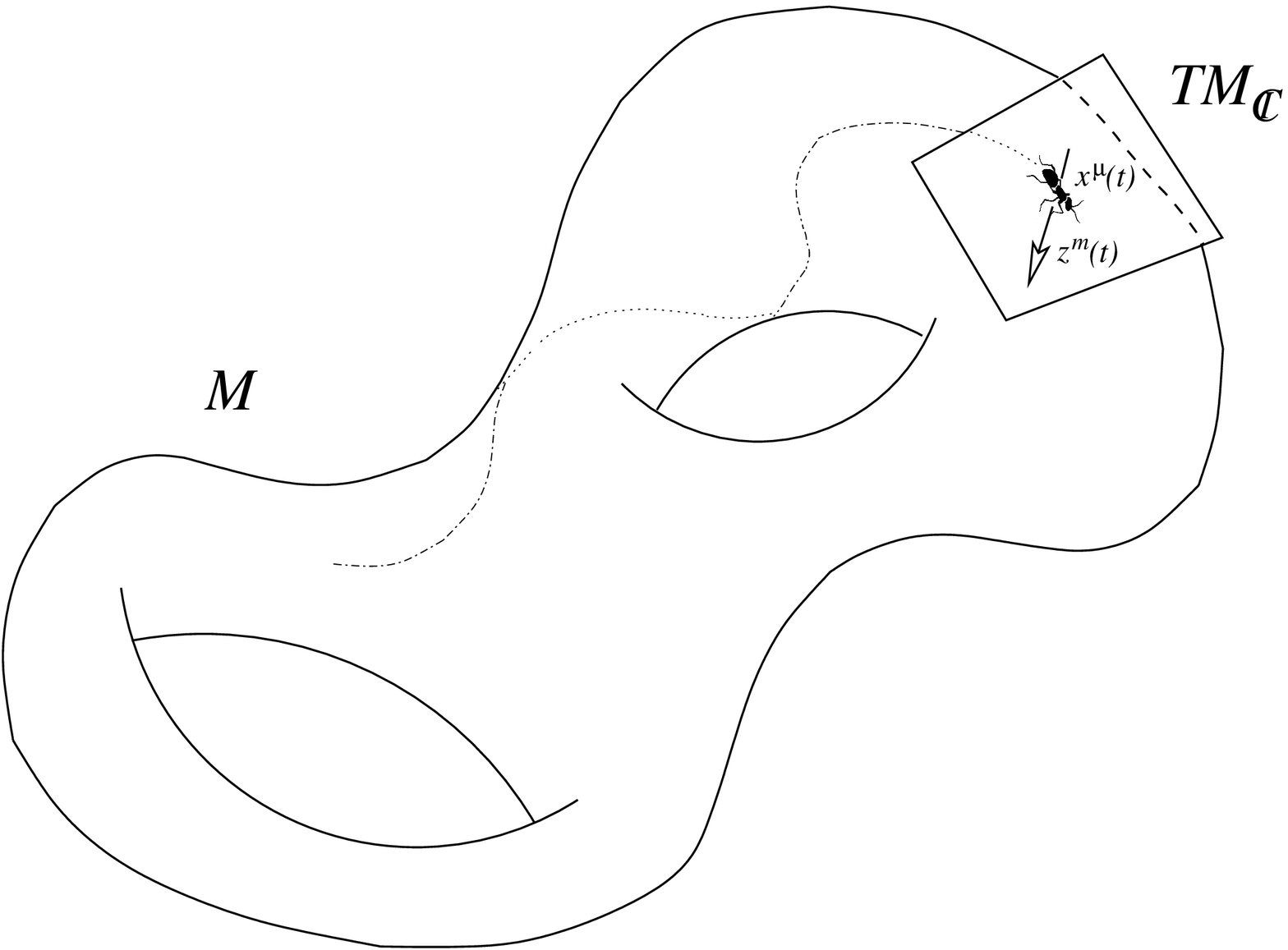}}
\caption{An ant laden with a complex tangent vector.\label{ant}}
\end{figure}

Let $(M,g_{\mu\nu})$ be an $n$-dimensional (pseudo-)Riemannian manifold with an orthonormal
frame $e^m$ so that\footnote{Although the metric signature
impacts the unitarity of the quantum Hilbert space of our model, all the results
presented here hold for arbitrary signature. Similarly, none of our results depend on the existence
of a global orthonormal frame.}
$ds^2=dx^\mu g_{\mu\nu}dx^\nu=e^m \eta_{mn} e^n$.
We consider the motion of an ant $x^\mu(t)$ -- as depicted in Fig.~\ref{ant} -- who carries a complex vector $z^m(t)$ (expressed relative to the orthonormal frame)
tangent to $M$. (In physics nomenclature, $z^m$ is referred to as commuting spinning degrees of freedom.) The Levi-Civita connection will be denoted by $\nabla$.
The ant determines its path and in which direction to hold the complex vector
by the system of generalized geodesic ODEs
\begin{gather}
\frac{\nabla\dot x^\mu}{dt}=\dot x^\nu R^\mu{}_\nu {}^m{}_n z^*_m z^n +
\nabla^\mu\! R^m{}_n{}^r{}_s\,  z^*_m z^n z^*_r z^s,\nonumber\\
\frac{\nabla z^m}{dt}=iR^m{}_n{}^r{}_s\, z^n z_r^* z^s .\label{evolve}
\end{gather}
The non-linear couplings to the curvature tensor on the right hand side of these
equations have been carefully chosen to maximize the set of constants of the motion.
They may obtained by extremizing a generalized energy integral
\begin{gather}
S=\int dt  \left\{ \frac12\dot x^\mu g_{\mu\nu} \dot x^\nu + i z_m^* \frac{\nabla z^m}{dt}+
\frac12 R^m{}_n{}^r{}_s \, z^*_m z^n z^*_r z^s\right\} .\label{action}
\end{gather}
To study constants of the motion, we look for symmetries of this action principle. The most obvious
of these are translations of $t\rightarrow t+\xi$ along the parameterized path traversed by our ant.
Inf\/initesimally this yields the invariance
\begin{gather*}
\delta x^\mu = \xi \dot x^\mu ,\qquad \delta z^m = \xi \dot z^m .
\end{gather*}
Less trivial, are symplectic transformations of $(z^m,z_m^*)$,
\begin{gather*}
\delta z^m  =-\alpha z^m +\beta z^{*m} ,\qquad
\delta z^{*m} = \gamma z^m + \alpha z^{*m} .
\end{gather*}
The parameters $(\alpha,\beta,\gamma)$ are real and correspond to the Lie algebra~$sp(2,{\mathbb R})$. The astute reader will observe that the symmetry transformation of $z^{*m}$ is
{\it not} the complex conjugate of $z^m$. Nonetheless, treating $z^m$ and $z^{*m}$ as
independent variables, the above $sp(2,{\mathbb R})$ transformations do leave the action invariant.
This is in fact suf\/f\/icient to ensure existence of corresponding constants of the motion and Noether charges. In the quantum theory, these charges will play an important geometric {\it r\^ole}.

The most interesting symmetries of the model interchange the complex tangent vector $z^m$
with the tangent vector $\dot x^\mu$ to the ant's path
\begin{gather}
\delta x^\mu = i (z^{*\mu} \varepsilon-z^{\mu} \varepsilon^*),\qquad
{\cal D} z^\mu = \dot x^\mu \varepsilon . \label{susy}
\end{gather}
Here, ${\cal D}$ is the covariant variation and is def\/ined by ${\cal D} v^\mu \equiv \delta v^\mu
+ \Gamma^\mu_{\rho\sigma} \delta x^\rho v^\sigma$ where $\Gamma$ denotes the Christof\/fel symbols.
It saves one from having to vary covariantly constant quantities.

The transformations~\eqn{susy} are not an exact symmetry for an arbitrary Riemannian
manifold. In fact, the action~\eqn{action} is invariant only when the
locally symmetric space condition
\begin{gather}
\nabla_\kappa R_{\mu\nu\rho\sigma}=0 ,\label{symmetric}
\end{gather}
holds. Or in other words, the Riemann tensor is covariantly constant. Constant curvature spaces provide
an, but by no means the only, example of such a manifold.

To compute constants of the motion we work in a f\/irst order formulation
\begin{gather*}
\dot x^\mu = \pi^\mu.
\end{gather*}
This and the evolution equations~\eqn{evolve} also follow from an
action principle
\begin{gather}
S^{(1)}=\int dt\left\{
p_\mu \dot x^\mu + i z^*_m \dot z^m -\frac12 \pi_\mu g^{\mu\nu} \pi_\nu
+\frac 12 R^m{}_n{}^r{}_s \, z^*_m z^n z^*_r z^s
\right\} ,\label{action1}
\end{gather}
where the covariant and canonical momenta $\pi_\mu$ and $p_\mu$ are related by
\begin{gather*}
\pi_\mu = p_\mu - i \omega_{\mu}{}^m{}_n z_m^* z^n .
\end{gather*}
Here the spin connection is determined by requiring covariant constancy of the
orthonormal frame $\nabla_\rho e_\mu{}^m= \partial_\rho e_\mu{}^m -\Gamma_{\rho\mu}^\nu e_\nu{}^m + \omega_\rho{}^m{}_n e_\mu{}^n=0$.

From the f\/irst order action~\eqn{action1} we immediately read of\/f the contact one-form
$p_\mu dx^\mu + i z_m^* dz^m$ which is already in Darboux coordinates, so
Poisson brackets follow immediately
\begin{gather*}
\{p_\mu,x^\nu\}_{\rm PB}=\delta_\mu^\nu ,\qquad
\{z^m,z^*_n\}_{\rm PB}=i\delta_n^m .
\end{gather*}
The Noether charges for the symmetries of the model can now be computed
\begin{gather}
H=\frac 12 \pi_\mu g^{\mu\nu} \pi_\nu -\frac 12 R^m{}_n{}^r{}_s \, z^*_m z^n z^*_r z^s ,\nonumber\\
f=\left(\begin{array}{cc} z^*_m z^{*m}&z^*_m z^m  \vspace{1mm}\\ z^*_m z^m& z_m z^m \end{array}\right) ,
\qquad
v=\ \left(\begin{array}{c}iz^{*\mu} \pi_\mu  \vspace{1mm}\\ i z^\mu \pi_\mu \end{array}\right) .
\label{Noether}
\end{gather}
The f\/irst of these is  the Hamiltonian. We have arranged the symplectic symmetry charges
in a~symmetric matrix~$f$ using the isomorphism between the symplectic Lie algebra and symmetric matrices. The remaining charges appear as a column vector~$v$ since they in fact form
a doublet representation of $sp(2,{\mathbb R})$.
It is important to remember that this latter pair of charges
are constants of the motion for locally symmetric spaces only.

\section{Quantization and geometry}
\label{quantum}

Quantization proceeds along usual lines replacing the Poisson brackets
by quantum commutators $[p,x]=-i\hbar$ and $[z,z^*]=\hbar$. We set
$\hbar=1$ in what follows and represent the canonical momentum
as a derivative acting on wavefunctions $\psi(x^\mu)$
\begin{gather*}
p_\mu=\frac 1i \frac{\partial}{\partial x^\mu} .
\end{gather*}
The spinning degrees of freedom become oscillators acting on a Fock space.
Rather than using the standard notation $z^m=a^m$ annihilating a Fock vacuum $a^m|0\rangle$,
we represent $|0\rangle=1$ and to preempt their geometric interpretation, set
\begin{gather*}
z^{*\mu}=dx^\mu ,\qquad z_{\mu}=\frac{\partial}{\partial(dx^\mu)} .
\end{gather*}
Therefore, wavefunctions become
\begin{gather*}
\Psi = \sum_{s=0}^{\infty} \psi_{\mu_1\ldots \mu_s}(x) dx^{\mu_1}\cdots dx^{\mu_s} ,
\end{gather*}
or in words -- sections of the symmetric tensor bundle~$SM$ over~$M$. Therefore, we can now start relating
quantum mechanical operations to dif\/ferential geometry ones on symmetric tensors.
Firstly, the quantum mechanical inner product yields the natural inner product for symmetric tensors
\begin{gather*}
\langle \Phi|\Psi\rangle = \int_M \sum_{s=0}^\infty s! \sqrt{g}\, \phi^{\mu_1\ldots\mu_s}\psi_{\mu_1\ldots\mu_s} .
\end{gather*}
Furthermore, the
covariant momentum corresponds to the covariant derivative $i\pi_\mu\Psi=\nabla_\mu \Psi$
(it is necessary to contract the open index $\mu$ with $dx^\mu$ for this to hold true for
subsequent applications of $\pi_\mu$).

Next we turn to the symplectic symmetries $f$ in~\eqn{Noether}. We call the of\/f-diagonal
charge
\begin{gather*}
{\bf N} = dx^\mu\frac{\partial}{\partial(dx^\mu)} ,
\end{gather*}
which simply counts the number of indices of a symmetric tensor
\begin{gather*}
{\bf N}\; \psi_{\mu_1\ldots \mu_s} dx^{\mu_1}\cdots dx^{\mu_s} = s \, \psi_{\mu_1\ldots\mu_s} dx^{\mu_1}\cdots dx^{\mu_s} .
\end{gather*}
We call the diagonal charges
\begin{gather*}
{\bf g}=dx^\mu g_{\mu\nu} dx^\nu ,\qquad
{\bf tr}=\frac{\partial}{\partial(dx^\mu)} \, g^{\mu\nu} \frac{\partial}{\partial(dx^\nu)} ,
\end{gather*}
as they produce new symmetric tensors by either multiplying
by the metric tensor and symmetrizing, or tracing a pair of indices
\begin{gather*}
{\bf g} \ \psi_{\mu_1\ldots \mu_s} dx^{\mu_1}\cdots dx^{\mu_s}
= g_{(\mu_1\mu_2}\psi_{\mu_3\ldots \mu_{s+2})} dx^{\mu_1}\cdots dx^{\mu_{s+2}},\\
{\bf tr} \ \psi_{\mu_1\ldots \mu_s} dx^{\mu_1}\cdots dx^{\mu_s}
=s(s-1) \psi^\mu{}_{\mu\mu_1 \ldots \mu_{s-2}}dx^{\mu_1}\cdots dx^{\mu_{s-2}} .
\end{gather*}
These three operators obey the $sp(2,{\mathbb R})$ Lie algebra
\begin{gather*}
[{\bf N},{\bf tr}]=-2\,{\bf tr} ,\qquad [{\bf N},{\bf g}]=2\,{\bf g} ,
\qquad
[{\bf tr},{\bf g}]=4\,{\bf N} + 2n .
\end{gather*}
We call its quadratic Casimir
\begin{gather*}
{\bf c}={\bf g}\, {\bf tr} -{\bf N}({\bf N}+n-2) .
\end{gather*}
Trace-free symmetric tensors with a def\/inite number of indices, $({\bf N}-s)\Psi=0={\bf tr}\Psi$,
are the highest weight vectors for unitary discrete series representations of this $sp(2,{\mathbb R})$
algebra.

The Noether charges $v$ in~\eqn{Noether} are linear in momenta and therefore covariant derivatives,
when acting on wavefunctions. We call them the gradient and divergence,
\begin{gather*}
{\bf grad}=dx^{\mu} \nabla_{\mu} ,\qquad
{\bf div}= \frac{\partial}{\partial(dx^\mu)} \nabla^\mu ,
\end{gather*}
because they are natural generalizations to symmetric tensors
of the exterior derivative and codif\/ferential for dif\/ferential forms.
To be sure
\begin{gather*}
{\bf grad} \, \psi_{\mu_1\ldots \mu_s} dx^{\mu_1}\cdots dx^{\mu_s}
= \nabla_{(\mu_1}\psi_{\mu_2\ldots \mu_{s+1})} dx^{\mu_1}\cdots dx^{\mu_{s+1}} ,\nonumber\\
{\bf div} \, \psi_{\mu_1\ldots \mu_s} dx^{\mu_1}\cdots dx^{\mu_s}
=s \nabla^\mu\psi_{\mu\mu_1 \ldots \mu_{s-1}}dx^{\mu_1}\cdots dx^{\mu_{s-1}} .
\end{gather*}
This pair of operators forms the def\/ining representation of $sp(2,{\mathbb R})$
\begin{gather*}
[{\bf N},{\bf grad}]={\bf grad} ,\qquad [{\bf N},{\bf div}]=-{\bf div} ,
\qquad
[{\bf tr},{\bf grad}]=2\, {\bf div} ,\qquad [{\bf div},{\bf g}]=2\, {\bf grad} .
\end{gather*}
It remains to commute the operators ${\bf div}$ and ${\bf grad}$.
The result is
\begin{gather}
[{\bf div},{\bf grad}]=\Delta - R^{\#\#} ,\label{divgrad}
\end{gather}
where $\Delta=\nabla^\mu\nabla_\mu$ is the Bochner Laplacian and
\begin{gather*}
R^{\#\#}\equiv R_{\mu}{}^\nu{}_\rho{}^{\sigma}dx^\mu\frac{\partial}{\partial(dx^\nu)}
dx^\rho\frac{\partial}{\partial(dx^\sigma)} .
\end{gather*}
This relation is closely analogous to that for the exterior derivative and codif\/ferential
$\{d,\delta\} = \Delta_F$ where $\Delta_F$ is the form-Laplacian. Here, since we are dealing
with symmetric tensors, the anticommutator is replaced by a commutator. A shrewd reader might
sense that the second order operator on the right hand side of~\eqn{divgrad} should be
related to the quantum mechanical Hamiltonian operator. This is indeed the case; calling
$\square = \Delta+R^{\#\#}$, we have
\begin{gather*}
[{\bf div},{\bf grad}]=\square - 2R^{\#\#} ,
\end{gather*}
where $\square = - 2 H$ so long as an appropriate operator ordering is chosen for the
Hamiltonian (a~full account is given in~\cite{Hallowell:2007qk}). On any manifold
\begin{gather*}
[\square,{\bf g}] = [\square,{\bf N}] = [\square,{\bf tr}]=0 .
\end{gather*}
Moreover, whenever the symmetric space condition~\eqn{symmetric}
holds, the operator~$\square$ is central
\begin{gather*}
[\square,{\bf div}]=[\square,{\bf grad}]=0 .
\end{gather*}
In fact, $\square$ is precisely the wave operator introduced quite some time ago
by Lichnerowicz on the basis of its special algebra with gradient and divergence operators~\cite{Lichnerowicz:1961}.
Finally, in the special case of constant curvature manifolds, choosing units in which
the scalar curvature $R=-n(n-1)$, the curvature operator $R^{\#\#}$ equals
the $sp(2,{\mathbb R})$ Casimir so that
\begin{gather}
[{\bf div},{\bf grad}]=\square - 2{\bf c} .\label{divgradc}
\end{gather}
If we include a further operator ${\bf ord}$ whose {\it r\^ole} is to count derivatives
\begin{gather*}
[{\bf ord},{\bf g}]=[{\bf ord},{\bf N}]=[{\bf ord},{\bf tr}]=0 ,
\\
[{\bf ord},{\bf grad}]={\bf grad} ,\qquad [{\bf ord},{\bf div}]={\bf div} ,
\qquad
[{\bf ord},\square]=2\square ,
\end{gather*}
then, $\{{\bf ord},{\bf g},{\bf N},{\bf tr},{\bf grad},{\bf div},\square\}$ form a maximal
parabolic subgroup of $sp(4,{\mathbb R})$ up to the rank~2 deformation by
the $sp(2,{\mathbb R})$ Casimir in~\eqn{divgradc}.
On f\/lat manifolds $M$, the operators $\{{\bf g},{\bf N},{\bf tr}$, ${\bf grad},{\bf div},\square\}$
obey the Fourier--Jacobi Lie algebra of $sp(2,{\mathbb R})$.
In the next Section, we present a novel reformulation of its universal enveloping algebra
based on introducing a certain square root of the Casimir operator ${\bf c}$.

\section[The Fourier-Jacobi algebra]{The Fourier--Jacobi algebra}
\label{FJalgebra}

\begin{figure}
\centerline{\includegraphics[height=8.5cm]{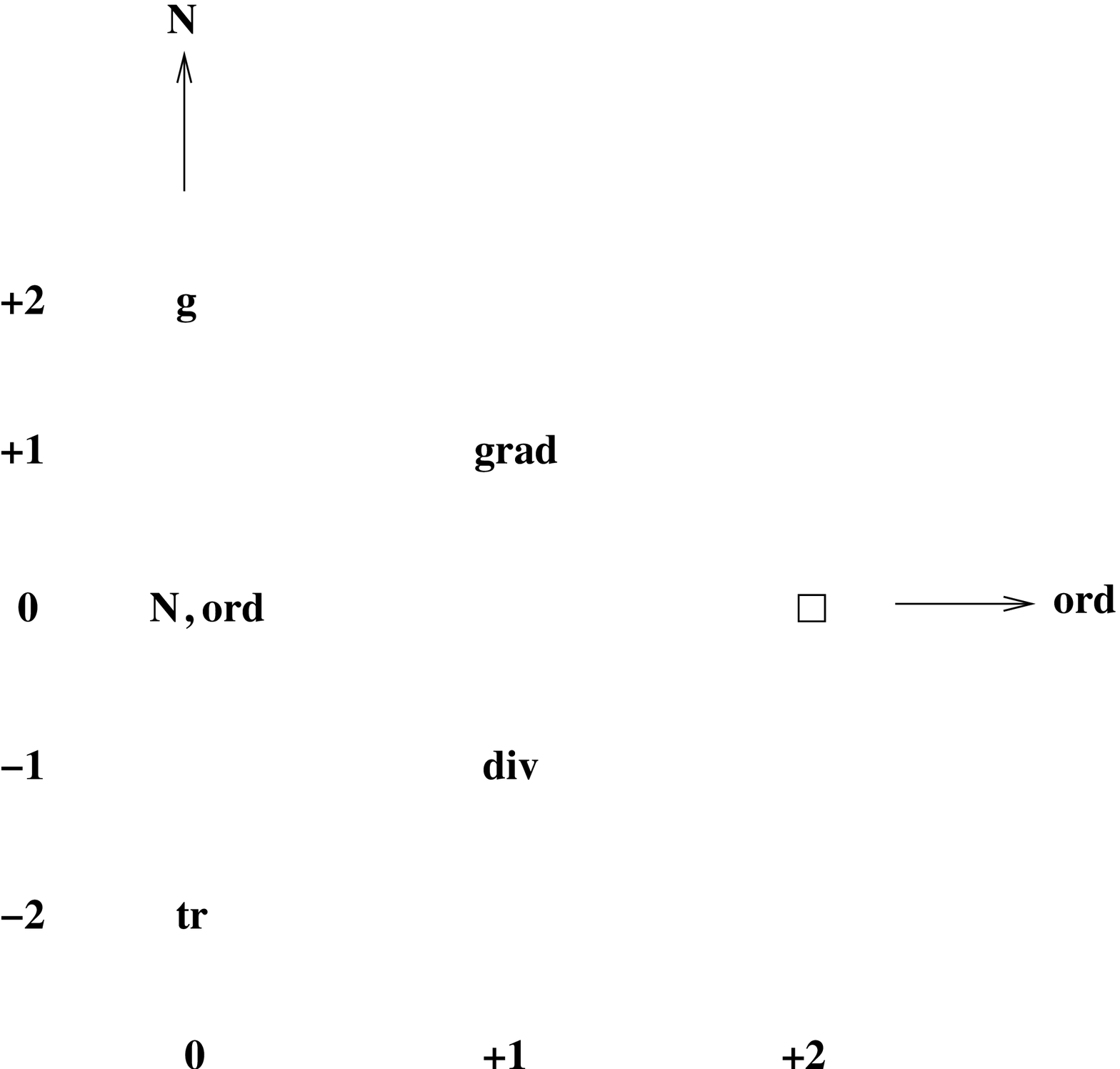}}
\caption{The root diagram for the Fourier--Jacobi Lie algebra.\label{roots}}
\end{figure}

Let us f\/irst collect together the deformed Fourier--Jacobi Lie algebra
built from geometric operators on constant curvature spaces
\begin{center}
\shabox{
\begin{tabular}{c}
$
[{\bf tr},{\bf g}]=4\,{\bf N} + 2n ,
$
\\[2mm]
$
\!\!\!\!\!\![{\bf N},{\bf tr}]=-2\,{\bf tr} ,\qquad\quad\ \ \   [{\bf N},{\bf g}]=2\,{\bf g} ,
$
\\[2mm]
$
[{\bf N},{\bf grad}]={\bf grad} ,\qquad\qquad [{\bf N},{\bf div}]=-{\bf div} ,
$
\\[2mm]
$
[{\bf tr},{\bf grad}]=2\, {\bf div} ,\qquad\qquad [{\bf div},{\bf g}]=2\, {\bf grad},
$
\\[2mm]
$
[{\bf div},{\bf grad}]=\square - 2{\bf c} .
$
\end{tabular}
}
\end{center}
Its root diagram is given in Fig.~\ref{roots}. From now on we compute
in the explicit realization given by its action on sections of the symmetric tensor bundle.
Therefore we are working with linear operators, so any algebra we f\/ind is automatically
consistent and associative.

To start with we analyze the $sp(2,{\mathbb R})$ Lie algebra built from $\{\bf g, N, tr\}$.
Unitary discrete series representations with respect to the adjoint involution
${\bf tr}^\dagger = {\bf g}$, ${\bf N^\dagger=N}$ are built from highest weights $\Phi$
such that
\begin{gather*}
{\bf N}\, \Phi=s\Phi ,\qquad {\bf tr}\, \Phi = 0 .
\end{gather*}
The highest weight module is spanned by
\begin{gather*}
\{\Phi, {\bf g}\, \Phi, {\bf g}^2 \Phi,\ldots\} .
\end{gather*}
We can characterize this representation by the eigenvalue $s$ of ${\bf N}$ acting on the
highest weight, or alternatively by the eigenvalue $-s(s+n-2)$ of the Casimir ${\bf c}$ acting
on any state in the module.

Conversely, given a eigenstate $\Psi$ of ${\bf N}$ and ${\bf c}$, we can determine which discrete series representation it belongs to by repeatedly applying the trace operator
\begin{gather*}
{\bf tr}^{k} \Psi \neq 0= {\bf tr}^{k+1} \Psi ,
\end{gather*}
which implies that $\Psi$ can be expressed in terms of a  highest weight vector
as~$\Psi={\bf g}^k \, \Phi$. {\it Our key observation is that it is highly advantageous
to introduce the linear operator ${\bm \kappa}$ whose eigenvalue acting on $\Psi$ is the depth $k$, namely}
\begin{gather}
{\bm \kappa} \equiv \frac{{\cal N - C}-1}2 ,\label{newop}
\end{gather}
where
\begin{gather*}
{\cal N}\equiv{\bf N}+\frac{n}{2} ,\qquad
{\cal C}\equiv\sqrt{\left(\frac{n-2}{2}\right)^2-{\bf c}} .
\end{gather*}
In other words, ${\bm \kappa}$ measures how far the symmetric tensor $\Psi$ is from being trace-free.
The operator $4\, {\cal N}$ is simply the right hand side of the ${\bf tr}$, ${\bf g}$ commutator.
More important is the square root of the Casimir ${\cal C}$ which acts on the highest weight $\Phi$
as
\begin{gather*}
{\cal C}\, \Phi = \left[s+\frac{n-2}2\right]\Phi ,
\end{gather*}
which explains equation~\eqn{newop}.

Our claim is that for many applications involving high powers of the operators ${\bf g}$, ${\bf N}$ and~${\bf tr}$,
rather than normal or Weyl orderings it is far more expeditious to work with functions of~${\cal N}$,~${\cal C}$
(up to perhaps an overall power of ${\bf g}$ or ${\bf tr}$). By way of translation, we note
that a normal ordered product of ${\bf g}$ and ${\bf tr}$ can be expressed as Pochhammer\footnote{Recall that the Pochhammer symbol is def\/ined as $(x)_m\equiv x(x+1)\cdots(x+m-1)$.}
functions of $({\cal N,C})$\ :
\begin{gather*}
{\bf g}\, {\bf tr}  =  ({\cal N -C} -1)({\cal N+C}-1) ,\qquad \Rightarrow \\
:({\bf g\, tr})^m: \equiv {\bf g}^m{\bf tr}^m\nn  = \frac1{4^m}\
\left(\frac{{\cal C-N}+1}2\right)_m
\left(\frac{{\cal C+N}-1}2\right)_m .
\end{gather*}
This claim has little signif\/icance until we introduce the doublet $({\bf grad,div})$.
Indeed, since the Casimir ${\bf c}$ has a rather unpleasant commutation
relation with either of these operators, computing in terms of $({\cal N,C})$
may seem unwise. In fact this is not the case once one appropriately modif\/ies the
divergence and gradient operators.

To motivate the claim we return to symmetric tensors. Suppose $\varphi_{\mu_1\ldots \mu_s}$
is trace-free, then its gradient $\nabla_{(\mu_1}\varphi_{\mu_2\ldots \mu_{s+1})}$ is in general not
trace free (unless the divergence of $\varphi$ happens to vanish). Since we would like to work with states
diagonalizing both ${\bf c}$ and ${\bf N}$, it is propitious to replace the regular gradient with its
trace-free counterpart
\begin{gather*}
\nabla_{(\mu_1}\varphi_{\mu_2\ldots\mu_{s+1})}-\frac{s}{2s+n-2}\ g_{(\mu_1\mu_2}
\nabla^\mu\varphi_{\mu_3\ldots\mu_{s+1})\mu} .
\end{gather*}
We denote this operator by $\wt{\bf grad}$. Having introduced ${\cal C}$ and ${\cal N}$, it has the simple
expression
\begin{gather}
\wt{\bf grad}   \equiv  {\bf grad} - {\bf g} \, {\bf div} \, \frac{1}{{\cal N+C}-1} ,\label{gradt}
\end{gather}
which we take to be its def\/inition acting on any section of the symmetric tensor bundle.
It is important to note that although this operator maps trace-free tensors
to trace-free tensors, it is designed to maintain how far a more general tensor
is from being trace-free. Therefore it does not project arbitrary tensors to trace-free ones.
We also introduce a similar def\/inition for a trace-free divergence following from
the quantum mechanical adjoint ${\bf grad}^\dagger = -{\bf div}$
\begin{gather}
\wt{\bf div}  \equiv\  {\bf div} - \frac{1}{{\cal N+C}-1} \, {\bf grad} \, {\bf tr}  .\label{divt}
\end{gather}
Note also that the linear operator ${{\cal N+C}-1}$ is indeed invertible since its spectrum is
$2s+2k+n-2$ on eigenstates ${\bf g}^k\Phi$
(these expressions also make sense in
dimensions
$n=1,2$ thanks to the operators ${\bf div}$ and ${\bf grad} \, {\bf tr}$).

The beauty of the operators $({\bf \wt{div},\wt{grad}})$ is that they {\it commute} with the
depth operator ${\bm \kappa}$. This implies
\begin{alignat*}{3}
&{\cal N} \, \wt{\bf div}=\wt{\bf div}\, ({\cal N}-1) ,\qquad && {\cal N}\, \wt{\bf grad}=\wt{\bf grad}\, ({\cal N}+1) ,&
\\
&{\cal C} \, \wt{\bf div} = \wt{\bf div} \, ({\cal C}-1) , \qquad && {\cal C}\, \wt{\bf grad}=\wt{\bf grad}\, ({\cal C}+1).&
\end{alignat*}
Moreover, an easy computation using the def\/initions~\eqn{gradt} and~\eqn{divt}
shows that the ordering of gradient and metric operators
can be interchanged at the cost of only a rational function of~$({\cal C, N})$
\begin{gather*}
\wt{\bf grad} \, {\bf g} = {\bf g}\, \wt{\bf grad} \, \frac{{\cal N+C}-1}{{\cal N+C}+1} ,
\qquad
{\bf tr} \, \wt{\bf div}  =  \frac{{\cal N+C}-1}{{\cal N+C}+1}\, \wt{\bf div} \, {\bf tr}   .
\end{gather*}
These relations allow us to invert equations~\eqn{gradt} and~\eqn{divt}
\begin{gather*}
{\bf grad}= \frac{1}{2}\, \frac{{\cal N+C}-3}{{\cal C}-1}\, \wt{\bf grad}
 + \frac1{2\cal C}\, \frac{{\cal N+C}-3}{{\cal N+C}-1}\, {\bf g}\, \wt{\bf div}
 ,\\
{\bf div} =
\frac12\, \wt{\bf div}\ \frac{{\cal N+C}-3}{{\cal C}-1}
+\frac12  \, \wt{\bf grad}\, {\bf tr}\, \frac{{\cal N+C}-3}{{\cal C(N+C}-1)} .
\end{gather*}
In turn we can now compute relations for reordering the gradient and trace operators
\begin{gather*}
{\bf tr}\, \wt{\bf grad} =
\frac{({\cal N+C}+1)({\cal N+C}-3)}{({\cal N+C}-1)^2}\, \wt{\bf grad}\ {\bf tr} ,\\
\wt{\bf div} \,  {\bf g} =
{\bf g} \, \wt{\bf div}\,
\frac{({\cal N+C}+1)({\cal N+C}-3)}{({\cal N+C}-1)^2} .
\end{gather*}
The f\/inal relation we need is for $\wt{\bf div}$ and $\wt{\bf grad}$. After some computations
we f\/ind
\begin{gather*}
\wt{\bf div}\, \wt{\bf grad} =
\frac{{\cal C}^2 ({\cal N+C}+1)({\cal N+C}-3)^2 }{({\cal C}+1)({\cal C}-1)({\cal N+C}-1)^3}\,
\wt{\bf grad}\, \wt{\bf div} \\
\phantom{\wt{\bf div}\, \wt{\bf grad} =}{}+
\frac{2 {\cal C}^2 ({\cal N+C}+1)}{({\cal C}+1)({\cal N+C}-1)^2}
\left[  \square + 2 \left({\cal C}+\frac n2 -1\right)\left({\cal C}-\frac n2 +1\right)\right] .
\end{gather*}
This result is valid in constant curvature spaces. The term in square brackets
equals $[\square - 2{\bf c}]$ and will be modif\/ied accordingly upon departure from constant
curvature.

\begin{figure}[t]
\begin{center}
\shabox{
\begin{tabular}{c}
\\[-5mm]
${\cal C}\, {\cal N}={\cal N}\, {\cal C}$\\[2mm]
${\bf tr}\, {\bf g} = ({\cal N+C}+1)({\cal N-C}+1)$\\[2mm]
${\bf g} \, {\bf tr} = ({\cal N+C}-1)({\cal N-C}-1)$\\[5mm]
$
{\cal N} \, {\bf tr}={\bf tr}\, ({\cal N}-2)\qquad\qquad {\cal N}\, {\bf g}={\bf g}\, ({\cal N}+2)
$\\[2mm]
$
 {\cal C}\, {\bf tr} = {\bf tr}\, {\cal C} \qquad\qquad {\cal C} \, {\bf g}={\bf g} \, {\cal C} $
\\[5mm]
\ \ \ \ $
{\cal N} \, \wt{\bf div}=\wt{\bf div}\, ({\cal N}-1)\qquad\qquad  {\cal N}\, \wt{\bf grad}=\wt{\bf grad}\, ({\cal N}+1)
$\\[2mm]
\ \ \ \ $
{\cal C} \, \wt{\bf div} = \wt{\bf div} \, ({\cal C}-1) \qquad \qquad {\cal C}\, \wt{\bf grad}=\wt{\bf grad}\, ({\cal C}+1)
$\\[5mm]
\ \ \ $\displaystyle
 \wt{\bf grad} \, {\bf g} = {\bf g}\ \wt{\bf grad} \, \frac{{\cal N+C}-1}{{\cal N+C}+1}$
\qquad\quad \
$\displaystyle
{\bf tr} \, \wt{\bf div}  =  \frac{{\cal N+C}-1}{{\cal N+C}+1}\, \wt{\bf div} \ {\bf tr}  \, \ \ \
$
\\[4mm]
$
\displaystyle  {\bf tr}\, \wt{\bf grad} =
\frac{({\cal N+C}+1)({\cal N+C}-3)}{({\cal N+C}-1)^2}\, \wt{\bf grad}\, {\bf tr}
$
\\[4mm]
$
\displaystyle \wt{\bf div} \,  {\bf g}  =
{\bf g} \, \wt{\bf div}\,
\frac{({\cal N+C}+1)({\cal N+C}-3)}{({\cal N+C}-1)^2}
$
\\[5mm]
$
\displaystyle \wt{\bf div}\ \wt{\bf grad}  =
\frac{{\cal C}^2 ({\cal N+C}+1)({\cal N+C}-3)^2 }{({\cal C}+1)({\cal C}-1)({\cal N+C}-1)^3}\
\wt{\bf grad}\ \wt{\bf div} \ \ \ \quad\qquad\qquad\quad\ \,
$
\\[4mm]
$
\displaystyle \phantom{\wt{\bf div}\ \wt{\bf grad}  =}{}
+
\frac{2 {\cal C}^2 ({\cal N+C}+1)}{({\cal C}+1)({\cal N+C}-1)^2}
\left[  \square + 2 \left({\cal C}+\frac n2 -1\right)\left({\cal C}-\frac n2 +1\right)\right]
$
\end{tabular}
}
\end{center}
\caption{Def\/ining relations for the $\wt{\, \cal U\, }\!(sp(2,{\mathbb R})^J)$ algebra.\label{thealgebra}}
\end{figure}

We  denote  the   new   algebra  built   from  $\{{\bf g,{\cal N},tr,{\cal C},\wt{grad},\wt{div},\square}\}$
 by     $\wt{\, \cal U\,}\!(sp(2,{\mathbb R})^J)$ and have collected together its def\/ining relations in Fig.~\ref{thealgebra}. As it is def\/ined by linear operators acting on symmetric tensors, associativity is
assured. An interesting, yet open, question is whether it can be def\/ined on the
universal enveloping algebra of $sp(2,{\mathbb R})^J$. Nonetheless, the explicit symmetric
tensor representation guarantees its consistency and therefore we may study it and its representations
as an abstract algebra in its own right.

\section{Conclusions}

We have presented a detailed study of symmetric tensors on curved manifolds.
The key technology employed is the quantum mechanics of a bosonic
spinning particle model. Also, many of our constructions were originally motivated
by studies of higher spin quantum f\/ield theories~\cite{Hallowell}. The spinning particle model presented here
is one of a general class of orthosymplectic spinning particle models that describe
spinors, dif\/ferential forms, multiforms, and indeed the most general tensor-spinor f\/ields on
a Riemannian manifold~\cite{Hallowell:2007qk}.

There are many applications and further research avenues. One simple question is that given the
strong analogy between the theory of dif\/ferential forms and the symmetric tensor one presented here,
are there symmetric tensor analogs of de Rham cohomology? The answer is yes. Recall, for example,
the Maxwell detour complex
\begin{gather}
\begin{array}{ccccccccccc}
  &                            &                        &           d              &
  &&&\delta&                         \\
0&\longrightarrow&\Lambda^0M&\longrightarrow&\Lambda^1M&\rightarrow\cdots \! \ \! \ \cdots
\rightarrow & \Lambda^{1}M&\longrightarrow& \Lambda^0M&\longrightarrow &0 .\\
&&&&\hspace{.07cm}|&\hspace{-1.22cm}\raisebox{.041cm}{\underline{\phantom{eeeseeyesyesyesyeseee}}}\hspace{-1.23cm}&\uparrow\hspace{.038cm} \\
&&&&&\delta d
\end{array}
\label{maxwell}
\end{gather}
This is mathematical shorthand for Maxwell's electromagnetism in curved backgrounds.
The physics translation is to replace the sequence of antisymmetric tensor bundles
$(\Lambda^0M,\, \Lambda^1M,\, \Lambda^1M$, $\Lambda^0M)$ by the words
\[
\mbox{\it (gauge parameters,
potentials, field equations, Bianchi identities)} .
\]
Then the fact that~\eqn{maxwell} is a complex implies that Maxwell's equations
$\delta d A=0$ are gauge invariant because $\delta d d=0$, and subject to a Bianchi identity
as $\delta\delta d=0$. An analogous  complex exists for symmetric tensors although not in general backgrounds, for brevity we give the f\/lat space result~\cite{Labastida:1987kw,Hallowell}
\begin{gather*}
\begin{array}{ccccccccccc}
  &                            &                        &           {\bf grad}              &
  &&&{\bf div}&                         \\
0&\longrightarrow&SM&\longrightarrow&SM&\rightarrow\cdots \ \! \ \cdots
\rightarrow & SM&\longrightarrow& SM&\longrightarrow &0 ,\\
&&&&\hspace{.07cm}|&\hspace{-1.22cm}\raisebox{.035cm}{\underline{\phantom{eeyeseyesyesyesyesee}}}\hspace{-1.24cm}&\uparrow\hspace{.038cm} \\
&&&&& G
\end{array}
\end{gather*}
where
\begin{gather*}
G=\square-\GRAD\, \DIV+\frac12\left(\GRAD^2\, \TR+\G\, \DIV^2\right)
    -\frac12\G\left(\square+\frac12\GRAD\, \DIV\right)\TR .
\end{gather*}
Notice that we have specif\/ied no grading on the symmetric tensor bundle~$SM$.
In fact the operator $G$ is the generating function for the equations of motion
(and actions) for massless higher spins of arbitrary degree.  A very fascinating question
is whether such complexes exist for the most general orthosymplectic spinning
particle models -- preliminary studies suggest an af\/f\/irmative answer~\cite{Deser}.

Another  interesting  open  question  is  the  generality  of  the  algebra
$\wt{\, \cal U\, }(sp(2,{\mathbb R})^J)$. For example, does there exist an algebra $ \wt{\, \cal U\, }(osp(2p|Q)^J)$
where $osp$ denotes the orthosymplectic superalgebra. Since the key idea is to include the
square root of the Casimir operator in the algebra, higher rank generalizations ought involve
the higher order Casimir operators. A positive answer to this question would be most
welcome and is under investigation~\cite{Cherney}.

\subsection*{Acknowledgements}
It is a pleasure to thank the organizers of the 2007 Midwest Geometry Conference, and especially
Susanne Branson for a truly excellent meeting in honor of Tom Branson.
We thank David Cherney, Stanley Deser, Rod Gover, Andrew Hodge, Greg Kuperberg,
Eric Rains and Abrar Shaukat for discussions.

\pdfbookmark[1]{References}{ref}
\LastPageEnding


\begin{thebibliography}{99}

\footnotesize\itemsep=0pt


\bibitem{Alvarez-Gaume:1983ig}
  Alvarez-Gaume L., Witten E., Gravitational anomalies,
  {\it Nuclear  Phys.~B} {\bf 234} (1984), 269--330.



\bibitem{Witten:1982im}
  Witten E., Supersymmetry and Morse theory,
  {\it J. Differential Geom.}  {\bf 17} (1982), 661--692.



\bibitem{Hallowell:2007qk}
  Hallowell K., Waldron A.,
  Supersymmetric quantum mechanics and super-Lichnerowicz algebras,
 {\it Comm. Math. Phys.},  to appear,
 \href{http://arxiv.org/abs/hep-th/0702033}{hep-th/0702033}.



\bibitem{Olver}
Olver P.J., Dif\/ferential hyperforms I, University of Minnesota, Report 82-101, 1982, 118 pages.\\
Olver P.J., Invariant theory and dif\/ferential equations, in Invariant Theory,
Editor S.~Koh,  {\it  Lecture Notes in Math.}, Vol.~1278, Springer-Verlag, Berlin, 1987, 62--80.



\bibitem{Henneaux}
  Dubois-Violette M., Henneaux M.,
  Generalized cohomology for irreducible tensor f\/ields of mixed Young symmetry type,
  {\it Lett.  Math. Phys.}  {\bf 49} (1999), 245--252,  \href{http://arxiv.org/abs/math.QA/9907135}{math.QA/9907135}.\\
Dubois-Violette M., Henneaux M.,  Tensor f\/ields of mixed Young symmetry type and $N$ complexes,
  {\it  Comm.  Math. Phys.}  {\bf 226} (2002), 393--418, \href{http://arxiv.org/abs/math.QA/0110088}{math.QA/0110088}.



\bibitem{Senovilla}
Edgar S.B., Senovilla J.M.M.,
A weighted de Rham operator acting on arbitrary tensor f\/ields and their local potentials,
 {\it J.  Geom. Phys.} {\bf 56} (2006), 2153--2162.



 \bibitem{Bekaert:2002dt}
 Bekaert X., Boulanger N.,
 Tensor gauge f\/ields in arbitrary representations of ${\rm GL}(D,{\mathbb R})$. Duality and Poincar\'e lemma,
 {\it Comm. Math. Phys.} {\bf 245}  (2004),  27--67,  \href{http://arxiv.org/abs/hep-th/0208058}{hep-th/0208058}.\\
Bekaert X., Boulanger N.,  On geometric equations and duality for free higher spins,
  {\it Phys. Lett.~B} {\bf 561} (2003), 183--190,   \href{http://arxiv.org/abs/hep-th/0301243}{hep-th/0301243}.\\
Bekaert X., Boulanger N.,   Tensor gauge f\/ields in arbitrary representations of ${\rm GL}(D,{\mathbb R})$,
 {\it Comm. Math. Phys.}  {\textbf{271}} (2007), 723--773,      \href{http://arxiv.org/abs/hep-th/0606198}{hep-th/0606198}.



\bibitem{deMedeiros:2003dc}
  de Medeiros P., Hull C.,
  Geometric second order f\/ield equations for general tensor gauge f\/ields,
   {\it JHEP} {\bf 2003} (2003), no.~5, 019, 26 pages, \href{http://arxiv.org/abs/hep-th/0303036}{hep-th/0303036}.\\
de Medeiros P., Hull C.,   Exotic tensor gauge theory and duality,
   {\it Comm.\ Math.\ Phys.}  {\bf 235} (2003), 255--273,   \href{http://arxiv.org/abs/hep-th/0208155}{hep-th/0208155}.



\bibitem{Lichnerowicz:1961}
Lichnerowicz A.,
Propagateurs et commutateurs en relativit\'e g\'en\'erale,
{\it Inst. Hautes \'Etudes Sci. Publ. Math.} (1961), no.~10, 56 pages.\\
Lichnerowicz A., Champs spinoriels et propagateurs en relativit\'e g\'en\'erale,
{\it Bull. Soc. Math. France} {\bf 92}  (1964), 11--100.



\bibitem{Hallowell}
Hallowell K., Waldron A.,
Constant curvature algebras and higher spin action generating functions,
 {\it Nuclear Phys.~B} {\bf 724} (2005), 453--486,  \href{http://arxiv.org/abs/hep-th/0505255}{hep-th/0505255}.



 \bibitem{Lecomte}
Duval C., Lecomte P., Ovsienko V.,
Conformally equivariant quantization: existence and uniqueness,
{\it Ann. Inst. Fourier (Grenoble)} {\bf 49} (1999), 1999--2029, \href{http://arxiv.org/abs/math.DG/9902032}{math.DG/9902032}.\\
Duval C.,  Ovsienko V.,
Conformally equivariant quantum Hamiltonians,
{\it Selecta Math. (N.S.)} {\bf 7} (2001), 291--320.



\bibitem{Labastida:1987kw}
   Labastida J.M.F.,
   Massless particles in arbitrary representations of the Lorentz group,
  {\it Nuclear Phys.~B} {\bf 322} (1989), 185--209.



\bibitem{Vasiliev:1988xc}
  Vasiliev M.A.,
  Equations of motion of interacting massless f\/ields of all spins as a free dif\/ferential algebra
  {\it Phys. Lett.~B} {\bf 209} (1988), 491--497.



\bibitem{Bekaert:2005vh}
  Bekaert X., Cnockaert S., Iazeolla C., Vasiliev M.A.,
  Nonlinear higher spin theories in various dimensions,
  \href{http://arxiv.org/abs/hep-th/0503128}{hep-th/0503128}.



\bibitem{Deser:2001pe}
  Deser S., Waldron A.,
 Gauge invariances and phases of massive higher spins in (A)dS,
  {\it Phys. Rev. Lett.}  {\bf 87} (2001), 031601, 4 pages, \href{http://arxiv.org/abs/hep-th/0102166}{hep-th/0102166}.\\
Deser S., Waldron A.,  Partial masslessness of higher spins in (A)dS,
  {\it Nuclear  Phys.~B} {\bf 607} (2001), 577--604,   \href{http://arxiv.org/abs/hep-th/0103198}{hep-th/0103198}.\\
Deser S., Waldron A.,   Stability of massive cosmological gravitons,
  {\it Phys.\ Lett.~B} {\bf 508} (2001), 347--353,   \mbox{\href{http://arxiv.org/abs/hep-th/0103255}{hep-th/0103255}}.\\
Deser S., Waldron A.,   Null propagation of partially massless higher spins in (A)dS and cosmological constant speculations,
  {\it Phys.\ Lett.~B} {\bf 513} (2001), 137--141,   \href{http://arxiv.org/abs/hep-th/0105181}{hep-th/0105181}.\\
Deser S., Waldron A.,    Arbitrary spin representations in de Sitter from dS/CFT with applications to dS supergravity,
  {\it Nuclear Phys.~B} {\bf 662} (2003), 379--392,  \href{http://arxiv.org/abs/hep-th/0301068}{hep-th/0301068}.\\
Deser S., Waldron A.,   Conformal invariance of partially massless higher spins,
  {\it Phys. Lett. B} {\bf 603} (2004), 30--34,  \href{http://arxiv.org/abs/hep-th/0408155}{hep-th/0408155}.



\bibitem{Deser}
Deser A., Waldron A., Einstein tensors for mixed symmetry higher spins, in preparation.

\bibitem{Cherney}
Cherney D., Hallowell D., Hodge A., Shaukat A., Waldron A., Higher rank Fourier--Jacobi algebras, in preparation.


\end{thebibliography}
\end{document}